\documentclass[11pt,reqno]{amsart}
\usepackage{amsmath,amssymb,verbatim,geometry,color,esint}
\usepackage[all]{xy}
\usepackage{mathrsfs}
\usepackage[pdftex,hyperindex]{hyperref}
\usepackage[pdftex]{graphicx}
\usepackage{tikz}
\usetikzlibrary{math}

\usepackage{accents}

\newcommand{\A}{\mathbb{A}}

\newcommand{\C}{\mathbb{C}}
\newcommand{\DD}{\mathbb{D}}
\newcommand{\G}{\mathbb{G}}

\renewcommand{\P}{\mathbb{P}}
\newcommand{\Q}{\mathbb{Q}}
\newcommand{\R}{\mathbb{R}}
\newcommand{\Z}{\mathbb{Z}}

\newcommand{\na}{\mathrm{na}}
\newcommand{\ra}{\mathrm{rad}}

\newcommand{\refe}{\mathrm{ref}}

\newcommand{\cE}{\mathcal{E}}
\newcommand{\cF}{\mathcal{F}}

\newcommand{\cH}{\mathcal{H}}

\newcommand{\cL}{\mathcal{L}}
\newcommand{\cM}{\mathcal{M}}
\newcommand{\cN}{\mathcal{N}}
\newcommand{\cO}{\mathcal{O}}
\newcommand{\cP}{\mathcal{P}}

\newcommand{\cS}{\mathcal{S}}

\newcommand{\cX}{\mathcal{X}}

\newcommand{\f}{\varphi}

\newcommand{\la}{\lambda}
\newcommand{\om}{\omega}
\newcommand{\Om}{\Omega}

\newcommand{\p}{\psi}

\newcommand{\hK}{\widehat{\mathrm{K}}}

\newcommand{{\EL}}{Euler--Lagrange}

\newcommand{\cf}{{\rm cf.\ }} 
\newcommand{\eg}{{\rm e.g.\ }} 
\newcommand{\ie}{{\rm i.e.\ }}

\newcommand{\be}{\textbf{e}}

\newcommand{\hto}{\hookrightarrow}

\newcommand{\hyb}{\mathrm{hyb}}

\newcommand{\dd}{\mathrm{d}}

\newcommand{\an}{\mathrm{an}}

\DeclareMathOperator{\ee}{E}
\DeclareMathOperator{\hh}{H}

\DeclareMathOperator{\mm}{M}
\DeclareMathOperator{\rr}{R}
\DeclareMathOperator{\pp}{P}

\DeclareMathOperator{\Aut}{Aut}

\DeclareMathOperator{\Ent}{Ent}

\DeclareMathOperator{\MA}{MA}

\DeclareMathOperator{\ord}{ord}

\DeclareMathOperator{\PSH}{PSH}
\DeclareMathOperator{\CPSH}{CPSH}

\DeclareMathOperator{\Ric}{Ric}

\DeclareMathOperator{\GL}{GL}

\DeclareMathOperator{\U}{U}

\DeclareMathOperator{\hcH}{\widehat{\mathcal{H}}}

\DeclareMathOperator{\Hnot}{H^0}

\DeclareMathOperator{\Lie}{\mathrm{Lie}}

\newcommand{\ddc}{dd^c}

\renewcommand{\div}{\mathrm{div}}

\newcommand{\triv}{\mathrm{triv}}

\newcommand{\nA}{non-Archimedean }
\newcommand{\ma}{Monge--Amp\`ere }

\newcommand{\simto}{\overset\sim\to}

\newcommand{\seb}{\color{blue}}

\numberwithin{equation}{section}       

\newtheorem{prop} {Proposition} [section]
\newtheorem{thm}[prop] {Theorem} 
\newtheorem{defi}[prop] {Definition}

\newtheorem{prop-def}[prop]{Proposition-Definition}

\newtheorem*{thmA}{Theorem A}

\theoremstyle{remark}

\title[Complex geometry and the YTD conjecture]{From complex to non-Archimedean geometry: an approach to the YTD conjecture}
\author{S\'ebastien Boucksom and Mattias Jonsson}
\date{\today}

\begin{document}
\begin{abstract}
    These notes expand on talks given by the authors at the 2025 Summer Research Institute in Algebraic Geometry in Fort Collins, Colorado. We discuss the relation between algebraic, analytic, and non-Archimedean geometry over the complex numbers, and sketch a proof of a version of the Yau--Tian--Donaldson conjecture for constant scalar curvature K\"ahler metrics.
\end{abstract}

\maketitle

\setcounter{tocdepth}{1}
\tableofcontents
%
%
%
%
\section*{Introduction}
Our objective in these notes is twofold. First, we want to explain some relations between three types of geometric objects defined over the complex numbers: complex algebraic varieties, complex manifolds, and certain Berkovich spaces.
Second, we will show how these relations can be used as a guiding principle towards a version of the Yau--Tian--Donaldson (YTD) conjecture, on an algebro-geometric  criterion for the existence of a constant scalar curvature K\"ahler metric in the first Chern class of an ample line bundle. 

\medskip
Consider a smooth projective complex variety $X$ together with an ample line bundle $L$, first viewed in the complex analytic category. The first Chern class $c_1(L)$ contains many K\"ahler forms, and a question going back to Calabi is whether we can find a `canonical' K\"ahler form $\omega\in c_1(L)$. For example, we can ask for (the K\"ahler metric corresponding to) $\omega$ to be a \emph{constant scalar curvature K\"ahler metric} (cscK metric for short). The cscK condition amounts to $\omega$ satisfying a nonlinear PDE; however, the \emph{YTD conjecture} asserts that the existence of a cscK metric in $c_1(L)$ is equivalent to an algebro-geometric condition on $(X,L)$. In the important special case when $X$ is Fano and $L=-K_X$, a cscK metric is the same as a \emph{K\"ahler--Einstein metric}, and the YTD conjecture was proved in~\cite{CDS15} (see also~\cite{Tia15}), with the algebro-geometric criterion given by \emph{K-polystability}. When the automorphism group $\Aut(X)$ is finite, K-polystability further simplifies to \emph{K-stability}. 

Here we will sketch a proof of the following result, a special case of~\cite[Theorem~A]{wytd}.

\begin{thmA}
    If $(X,L)$ is a polarized smooth projective complex variety with $\Aut^0(X,L)$ trivial, then there exists a cscK metric in $c_1(L)$ iff $(X,L)$ is $\hK$-stable.
\end{thmA}
Here $\Aut^0(X,L)$ denotes the identity component of the automorphism group of $(X,L)$, and triviality means that it reduces to the fiberwise scaling action of $\C^\times$. 

Here $\hK$-stability is a condition of algebro-geometric nature, but formulated in terms of \emph{non-Archimedean geometry} in the sense of Berkovich~\cite{BerkBook}. Namely, to the polarized variety $(X,L)$ we can associate another pair $(X_\na,L_\na)$, consisting of the Berkovich analytifications of $X$ and $L$ with respect to the (non-Archimedean) \emph{trivial absolute value} on $\C$, in which every non-zero complex number has norm 1.

The space $X_\na$ is a compact Hausdorff space, and comes with a structure sheaf (that we shall not directly make use of). It contains as a dense subset the set of \emph{divisorial valuations}, associated to prime divisors on birational models of $X$, and each element of $X_\na$ more generally corresponds to a valuation on the function field of a subvariety of $X$. A point emphasized by Berkovich himself is that spaces such as $X_\na$ and $L_\na$ have many features in common with their complex analytic counterparts. Our approach to the YTD conjecture very much uses this point of view. 

As opposed to the proof in~\cite{CDS15}, the approach in~\cite{YTD} and, more generally, the proof of Theorem~A, is variational in nature. Very roughly, the idea is to study the space $\cH$ of smooth positive metrics on $L$. Taking the curvature form yields a bijection between $\cH/\R$ and the space of K\"ahler forms in $c_1(L)$, and cscK metrics correspond to minimizers of a certain functional $\mm\colon\cH\to\R$, the \emph{Mabuchi functional}.

In order to reveal a key geodesic convexity property of the Mabuchi functional, one extends it to the completion of $\cH$ with respect to a `tautological' $L^1$-metric first studied by Darvas~\cite{Dar15}. This completion $\hcH$, and functionals upon it, can be studied using \emph{global pluripotential theory}; in fact, $\hcH$ can be identified with the space $\cE^1$ of plurisubharmonic (psh) metrics on $L$ of finite energy.

Very roughly speaking, $\hcH$ behaves like an infinite-dimensional symmetric space of non-compact type, in accordance with Donaldson's vision for $\cH$ in~\cite{Don99}. In particular, $\hcH$ is uniquely geodesic. Moreover, the Mabuchi functional extends canonically to $\hcH$ and is geodesically convex~\cite{BerBer,BDL17}. Further, minimizers of the Mabuchi functional on $\hcH$ must lie in $\cH$ and thus correspond to cscK metrics~\cite{CC2}.

Using these and other nontrivial results, one can show that a cscK metric in $c_1(L)$ exists iff $(X,L)$ is \emph{geodesically stable} in the sense that the Mabuchi functional grows along all nontrivial geodesic rays. Still, this is not a condition of algebro-geometric nature.

At this point it is useful to compare $\hcH$ with an actual symmetric space of non-compact type, namely the space $\cN(V)$ of Hermitian norms $\|\cdot\|=e^{-\chi}$ on a finite-dimensional vector space $V$. This space is uniquely geodesic, and the set of (constant-speed) geodesic rays $\{\chi_t\}_{t\ge0}$ emanating from a fixed norm $\chi_0$ can be identified with the space $\cN_\na(V)$ of \emph{non-Archimedean} norms on $V$, by mapping $\{\chi_t\}$ to $\chi\in\cN_\na(V)$ defined by $\chi(v):=\lim_{t\to\infty}t^{-1}\chi_t(v)$. 

Following Berkovich's philosophy, we developed in~\cite{trivval,nakstab1} global pluripotential theory over $\C$ equipped with the trivial valuation; the main results are remarkably similar to the ones in the complex analytic setting. A non-Archimedean version of $\cH$ is the space $\cH_\na$ of (non-Archimedean) \emph{Fubini--Study metrics} on $L_\na$, which furthermore stands in bijection with the space of (normal, ample) \emph{test configurations} for $(X,L)$. We can equip $\cH_\na$ with a Darvas metric, and the completion $\hcH_\na$ can be identified with the space $\cE^1_\na$ of non-Archimedean psh metrics of finite energy. As in the case of norms above, $\hcH_\na$ can be embedded into the space of geodesic rays in $\hcH$, but now we get a proper subset.

The Mabuchi functional has a trivially valued analog $\mm_\na\colon\cH_\na\to\R$, which can be used to define K-stability, as explained in~\cite{BHJ1}. This functional furthermore extends naturally to $\hcH_\na$, and we now say that $(X,L)$ is \emph{$\hK$-stable} if $\mm_\na\ge0$ on $\hcH_\na$, with strict inequality outside the set of constant potentials. 

The final piece of the proof of Theorem~A is now to prove that geodesic stability is equivalent to $\hK$-stability, and for this we need to compute the slope of the Mabuchi functional along a (psh) geodesic ray $\{\phi_t\}_{t\ge0}$ in $\hcH$ emanating from a fixed reference metric $\phi_0$. By a key result of C.~Li~\cite{LiGeod}, the convex function $\R_{\ge 0}\ni t\mapsto\mm(\phi_t)$ has infinite slope at infinity when $\{\phi_t\}$ lies outside $\cE^1_\na$. Otherwise, this ray defines an element $\f\in\cE^1_\na$; strengthening Li's results, we then prove that the slope at infinity satisfies 
\[
\lim_{t\to\infty} t^{-1}\mm(\phi_t)=\mm_\na(\f).
\]
By definition, $\hK$-stability requires us to verify a condition for all elements in the space $\hcH$, which is quite large. However, $\hK$-stability turns out (see Theorem~\ref{thm:ytd1}) to be equivalent to \emph{uniform} $\hK$-stability, which can be tested on filtrations of the section ring of $L$, or rational convex combinations of divisorial valuations on $X$.

This note is organized as follows. We start in~\S\ref{sec:norms} by briefly studying the space of Hermitian and non-Archimedean norms on a finite-dimensional complex vector space, following~\cite{BE}. Then, in~\S\ref{sec:archgeom} and~\S\ref{sec:nageom} we note some basic facts on the complex-analytic and Berkovich analytification of a  polarized smooth projective variety. Section~\ref{sec:gppt} is devoted to global pluripotential theory, introducing the spaces $\cH\subset\hcH\simeq\cE^1$ and $\cH_\na\subset\hcH_\na\simeq\cE^1_\na$, and functionals thereupon. In~\S\ref{sec:E1geom} we describe how $\cE^1$ has some features of a symmetric space, and how $\cE^1_\na$ sits inside the space of geodesic rays in $\cE^1$. We also review the crucial result on the asymptotics of the Mabuchi functional along geodesic rays. Finally, in~\S\ref{sec:ytd} we outline the proof of the YTD conjecture and in~\S\ref{sec:extensions} we make some final remarks.

\medskip\noindent\textbf{Acknowledgements}. These notes follow talks by the authors at the 2025 Summer Research Institute in Algebraic Geometry in Fort Collins, Colorado. We are grateful to the organizers for the opportunity to present our work. We also thank Ying Wang for comments, and the referee for useful comments and suggestions. The first author was partially supported by the ANR project AdAnAr (ANR-24-CE40-6184) and the KRIS project. The second author was partially supported by NSF grants DMS-2154380 and DMS-2452797. 

%
%
%
%
\section{From Hermitian to non-Archimedean norms}\label{sec:norms}
We start by describing a simple but illustrative situation where the Archimedean and non-Archimedean worlds meet: two spaces of norms on a finite-dimensional vector space. We refer to~\cite{ICM,BE} for details on the material here.

Fix a complex vector space $V$ of finite dimension $n\ge1$. 
%
\subsection{The space of Hermitian norms} 
Let $\cN(V)$ be the space of Hermitian norms $\|\cdot\|$ on $V$. We write such norms `additively', using $\chi:=-\log\|\cdot\|$. Given a basis $\be=(e_1,\dots,e_n)$ of $V$ and $\la\in\R^n$, we have a unique norm $\chi\in\cN(V)$ that is diagonalized in the basis $\be$ and satisfies $\chi(e_i)=\la_i$ for all $i$. This yields an injection
$$
\iota_\be\colon\R^n\hto\cN(V),
$$
whose image, called the \emph{flat} associated to $\be$, consists in all Hermitian norms that are diagonalized in the given basis. The terminology comes from the fact that $\cN(V)\simeq\GL(n,\C)/\U(n)$ can be viewed as a Riemannian symmetric space, whose maximal totally geodesic submanifolds isometric to a Euclidean space are precisely the above flats. 

Any element of $\cN(V)$ is contained in a flat. In fact, as any two Hermitian norms can be simultaneously diagonalized, any two elements of $\cN(V)$ belong to a common flat. As a consequence, the length metric $\dd_2$ of the Riemannian symmetric space $\cN(V)$ is characterized as the unique metric whose pullback under each $\iota_\be\colon\R^n\hto\cN(V)$ coincides with the (Euclidean) $\ell^2$-distance. 

\subsection{Finsler metrics}
More generally, consider a norm $\tau$ on $\R^n$ that is symmetric, \ie invariant under the action of $\cS_n$, such as the $\ell^p$-norm for $p\in [1,\infty]$. Then $\tau$ can be used to define a natural Finsler norm on the tangent bundle of the manifold $\cN(V)$, whose induced length (pseudo)metric $\dd_\tau$ is characterized as the unique metric on $\cN(V)$ that pulls back to $\tau$ under every embedding $\iota_\be\colon\R^n\hto\cN(V)$. 

When $\tau$ is the $\ell^2$-norm, we thus recover the Riemannian metric $\dd_2$ above. By equivalence of norms on $\R^n$, all $\dd_\tau$ metrics on $\cN(V)$ are  Lipschitz equivalent, and the metric space $(\cN(V),\dd_\tau)$ is thus complete and locally compact. When $\tau$ is the $\ell^\infty$-norm, $\dd_\tau=\dd_\infty$ is the norm given by $\dd_\infty(\chi,\chi')=\sup_{v\in V}|\chi(v)-\chi'(v)|$.

Every straight line in a flat is a $\dd_\tau$-geodesic, and we get in this way a distinguished class of geodesics with respect to which the metric space $(\cN(V),\dd_\tau)$ is uniquely geodesic. Note that other, non-distinguished $\dd_\tau$-geodesics also exist, when the norm $\tau$ on $\R^n$ is not strictly convex (\eg for the $\ell^1$-norm).  

The metric space $(\cN(V),\dd_\tau)$ is furthermore \emph{Busemann convex} with respect to distinguished geodesics: for any two such geodesics $\{\chi_t\}_{t\in I}$, $\{\chi'_t\}_{t\in I}$, $\dd_\tau(\chi_t,\chi'_t)$ is a convex function of $t\in I$. When $\tau$ is the $\ell^2$-norm, this boils down to the fact that the Riemannian symmetric space $\cN(V)\simeq\GL(n,\C)/\U(n)$ has nonpositive sectional curvature.

\subsection{The radial limit space}\label{sec:radial}
A distinguished geodesic ray $\{\chi_t\}_{t\ge 0}$ in $(\cN(V),\dd_\tau)$ is, by definition, a linear ray in some flat of $\cN(V)$, and hence of the form $\chi_t=\iota_\be(t\la)$ for some basis $\be$ and $\la\in\R^n$. 

Two such rays $\{\chi_t\}$, $\{\chi'_t\}$ are \emph{parallel} if $\dd_\tau(\chi_t,\chi'_t)=O(1)$, a condition that does not depend on the choice of $\tau$ by Lipschitz equivalence of the $\dd_\tau$-metrics. Parallelism is an equivalence relation on the set of geodesic rays, the equivalence class of a ray being called its \emph{direction}, and the \emph{radial limit space} $\cN_\ra(V)$ of $\cN(V)$ is defined as the set of directions of (distinguished, constant speed) geodesic rays $\{\chi_t\}_{t\ge 0}$ in $\cN(V)$. Given a base point $\chi_0\in\cN(V)$, any element of $\cN_\ra(V)$ is the direction of a unique distinguished geodesic ray emanating from $\chi_0$. 

The Busemann convexity of $\dd_\tau$ allows us to equip $\cN_\ra(V)$ with the \emph{radial metric} 
\[
\dd_{\tau,\ra}\left(\{\chi_t\},\{\chi'_t\}\right):=\lim_{t\to\infty}t^{-1}\dd_\tau(\chi_t,\chi'_t).
\]
The resulting metric space $(\cN_\ra(V),\dd_{\tau,\ra)}$ is then  complete, and also uniquely geodesic and Busemann convex with respect to a natural distinguished class of geodesics. Its Lipschitz equivalence class is again independent of $\tau$, but the underlying topological space is \emph{not} locally compact as soon as $n=\dim V>1$; see below.

We emphasize that our geodesics have constant speed not necessarily equal to $1$. Time reparametrization thus yields a natural scaling action of $\R_{>0}$ on $\cN_\ra(V)$, with respect to which the radial metric $\dd_{\tau,\ra}$ is homogeneous. In particular, $(\cN_\ra(V),\dd_{2,\ra})$ is isomorphic to the metric cone over the \emph{visual boundary} of the Riemannian symmetric space $\GL(n,\C)/\U(n)$, equipped with the \emph{Tits metric}. 

For completeness, we show that $(\cN_\ra(V),\dd_{\tau,\ra})$ is not locally compact when $n>1$. Using scaling, it suffices to prove that the closed unit ball in $(\cN_\ra(V),\dd_{\infty,\ra})$ with center at the constant geodesic $0\in\cN_\ra(V)$ fails to be compact. Pick any ordered basis $\be=(e_1,\dots,e_n)$ for $V$, and set $\be^{(m)}:=(e_1-me_2,e_2,e_3,\dots,e_n)$ for $m\ge1$. For $t\ge0$, set $\chi_t^{(m)}=\iota_{\be^{(m)}}(t,0,\dots,0)\in\cN(V)$. Then $\{\chi^{(m)}_t\}\in\cN_\ra(V)$, and $d_{\infty,\ra}(0,\{\chi^{(m)}_t\})=1$ for all $m$, where $0\in\cN_\ra(V)$ is the constant geodesic ray. Now, for any nonzero vector $v\in V$, there exists $m_0\ge1$ such that $\chi^{(m)}_t(v)=0$ for $m\ge m_0$ and all $t$. If some subsequence $\{\chi^{(m_j)}_t\}$ had a limit in $\cN_\ra(V)$, this limit would therefore equal 0, contradicting $\dd_{\infty,\ra}(\{\chi^{(m_j)}\},0)=1$.
%
%
\subsection{The space of non-Archimedean norms}\label{sec:nanorms}
The radial limit space $\cN_\ra(V)$ admits a nice description in terms of \emph{non-Archimedean norms} on $V$ with respect to the trivial valuation on $\C$. Adopting again `additive' notation, these are defined as functions $\chi\colon V\to\R\cup\{+\infty\}$ that satisfy for all $v,w\in V$
\begin{itemize}
    \item $\chi(v+w)\ge\min\{\chi(v),\chi(w)\}$; 
    \item $\chi(av)=\chi(v)$ for $a\in\C^\times$; 
    \item $\chi(v)=+\infty$ iff $v=0$.
\end{itemize}
Setting $\cF^\la V:=\{\chi\ge\la\}$ for $\la\in\R$ shows that a non-Archimedean norm $\chi$ is equivalently given by a decreasing \emph{filtration} $(\cF^\la V)_{\la\in\R}$ of $V$ by linear subspaces that is exhaustive and separating, \ie $\cF^\la V=V$ for $\la\ll0$ and $\cF^\la V=0$ for $\la\gg0$, and also left-continuous, \ie $\cF^\la V=\bigcap_{\mu<\la}\cF^\mu V$ for all $\la\in\R$.

Denote by $\cN_\na(V)$ the set of non-Archimedean norms. Given a basis $\be=(e_1,\dots,e_n)$ and $\la\in\R^n$ there exists a unique norm $\chi\in\cN_\na(V)$ which is diagonalized by $\be$ and such that $\chi(e_i)=\la_i$ for all $i$, \ie
$$
\chi(\sum_ia_ie_i)=\min\{\la_i\mid a_i\ne0\}
$$
for all $a_i\in\C$. This defines an injection $\iota_{\be,\na}\colon\R^n\hto\cN_\na(V)$, whose image is now called an \emph{apartment}. 

As in the Hermitian case, any two elements in $\cN_\na(V)$ belong to a common apartment. For any symmetric norm $\tau$ on $\R^n$ there exists a unique metric $\dd_{\tau,\na}$ on $\cN_\na(V)$ whose pullback under each $\iota_{\be,\na}\colon\R^n\to\cN_\na(V)$ coincides with $\tau$. As it turns out, the resulting metric space $(\cN_\na(V),\dd_{\tau,\na})$ is isometrically isomorphic to the radial limit space $(\cN_\ra(V),\dd_{\tau,\ra})$. 

To see this, pick a distinguished geodesic ray $\{\chi_t\}$, and write $\chi_t=\iota_\be(t\la)$ for some basis $\be=(e_1,\dots,e_n)$ for $V$ and some $\la\in\R^n$. Then the non-Archimedean norm $\chi:=\iota_{\be,\na}(\la)$ satisfies 
$$
\chi(v)=\lim_{t\to\infty} t^{-1}\chi_t(v)
$$
for any $v\in V$, and hence only depends on the direction of $\{\chi_t\}$. For any symmetric norm $\tau$ on $\R^n$, it is now easy to see that this construction gives rise to an isometric isomorphism
$$
(\cN_\ra(V),\dd_{\tau,\ra})\simto (\cN_\na(V),\dd_{\tau,\na}),
$$
equivariant under the scaling action of $\R_{>0}$. 
%
\subsection{Generalization}
The situation above can be generalized as follows. Let $G$ be a compact connected Lie group with complexification $G_\C$. The choice of a $G$-invariant scalar product on $\Lie G$ turns $\Sigma(G):=G_\C/G$ into a Riemannian symmetric space. 

Given a maximal compact torus $T\subset G$, $\Lie T\simeq\Sigma(T)$ embeds as a flat of $\Sigma(G)$, \ie a maximal submanifold isometric to Euclidean space. Translating by $g\in G_\C$ yields an embedding 
$$
\iota_g\colon\Lie T\hto\Sigma(G),
$$
onto a flat of $\Sigma(G)$, and all flats arise in this way. 

More generally, any norm $\tau$ on $\Lie T$ invariant under the (finite) Weyl group $W=N_G(T)/T$ induces a $G$-invariant norm on $\Lie G$, and hence a Finsler norm on the tangent bundle of $\Sigma(G)$. The associated length metric $\dd_\tau$ is characterized as the unique metric on $\Sigma(G)$ whose pullback under each $\iota_g$ coincides with $\tau$. 

Straight lines in $\Lie T$ give rise to a distinguished class of geodesics in $(\Sigma(G),\dd_\tau)$ with respect to which it is uniquely geodesic and Busemann convex. It thus admits a radial limit space $(\Sigma_\ra(G),\dd_{\tau,\ra})$, which coincides with the metric cone over the visual boundary. 

The analog of the space of non-Archimedean norms is the \emph{conical Tits building} $\Sigma_\na(G)$. By construction, it is a union of apartments, \ie images of injections $\iota_{g,\na}\colon\Lie T\hto\Sigma_\na(G)$ parametrized by $g\in G_\C$, and admits a unique metric $\dd_{\tau,\na}$ whose pullback under each $\iota_{g,\na}$ coincides with $\tau$. We then have a canonical $\R_{>0}$-equivariant isometric isomorphism 
$$
(\Sigma_\ra(G),\dd_{\tau,\ra})\simto(\Sigma_\na(G),\dd_{\tau,\na}). 
$$
We refer to~\cite[Appendix~A]{wytd} for the details, which are well-known at least when $\tau$ is a Euclidean norm.

%
%
%
%
\section{Complex algebraic and analytic geometry}\label{sec:archgeom}
Let $X$ be a smooth complex algebraic variety. We use additive conventions for line bundles, writing $L+M$ instead of $L\otimes M$. 
%
%
\subsection{Analytification}
The algebraic variety $X$ induces a complex manifold, its \emph{analytification}, that we slightly abusively also denote by $X$. Similarly, any line bundle $L$ on $X$ induces a holomorphic line bundle, also denoted $L$. At least when $X$ is proper, identifying $X$ with its analytification is quite harmless, by the GAGA principle. On the other hand, viewing $X$ as a complex manifold leads to new tools, such as Hodge theory.
%
%
%
%
\subsection{Metrics on line bundles}\label{sec:ametrics}
For a line bundle $L$ we denote by $L^\times$ the total space of $L$ with the zero section removed. By a \emph{metric} (resp.\ \emph{singular metric}) on $L$ we mean a function $\phi\colon L^\times\to\R$ (resp.\ $\phi\colon L^\times\to\R\cup\{\pm\infty\}$) such that $\phi(\la v)=\phi(v)-\log|\la|$ for $v\in L^\times$ and $\la\in\C^\times$. As $L^\times$ is a complex manifold, there are natural notions of continuous and smooth ($C^\infty$) metrics on $L$.

The space of metrics on $L$ admits various natural operations. For example, we can take the maximum (or minimum) of finitely many metrics, and add constants to metrics. If $\phi_i$ is a metric on $L_i$, $i=1,2$, then $\phi_1+\phi_2$ is a metric on $L_1+L_2$. If $\phi$ is a metric on $L$, then any morphism $f\colon Y\to X$ induces a metric $f^\star\phi$ on $f^\star L$. Metrics on $\cO_X$ can be viewed as functions on $X$, by evaluating on the canonical section 1. These remarks also apply to singular metrics.

Any global section $s$ of $L$ defines a singular metric $\log|s|$ on $L$, whose restriction to the fiber $L^\times(x)$ at $x\in X$ is given by $\log|s(x)/\cdot|$. This metric is smooth outside the zero locus of $s$. A metric $\phi$ on an ample line bundle $L$ is a \emph{Fubini--Study metric} if there exist $m\ge 1$ and global sections $s_1,\dots,s_N$ of $mL$ that define an embedding of $X$ into $\P^{N-1}$, such that 
\[
\phi=\frac1{2m}\log\sum_{i=1}^N|s_i|^2.
\]

The \emph{curvature form} $\ddc\phi$ of a smooth metric $\phi$ on $L$ is defined by $\ddc\phi=\frac{\sqrt{-1}}{\pi}\partial\overline\partial(\phi\circ s)$
for any local nonvanishing holomorphic section $s$ of $L$. This is a closed form in $c_1(L)$. We say that $\phi$ is \emph{positive} if $\ddc\phi$ is positive, that is, a K\"ahler form. For example, a Fubini--Study metric on an ample line bundle $L$ is smooth and positive. In fact, the existence of a smooth positive metric on $L$ is equivalent to $L$ being ample, by the Kodaira embedding theorem. 

When $L$ is ample, we denote by $\cH=\cH_L$ the space of positive smooth metrics on $L$. By the $\ddc$-lemma, the curvature form gives a bijection between $\cH/\R$ and the space of K\"ahler forms in $c_1(L)$.
%
%
\subsection{Canonical metrics}
Assume $X$ is projective and $L$ ample. Calabi asked whether the cohomology class $c_1(L)$ contains a `canonical' K\"ahler form $\omega$. For example, we may look for \emph{constant scalar curvature K\"ahler metrics} (cscK metrics), satisfying $\Ric\omega\wedge\omega^{n-1}=c\cdot\omega^n$, for an algebraically determined constant $c\in\R$, where $\Ric\omega\in c_1(X)$ is the Ricci form. More general classes of canonical metrics include \emph{extremal metrics}, but we shall largely stick to the cscK case in these notes. If the canonical bundle $K_X$ is proportional to $L$, $K_X=\lambda L$, then a cscK metric is in fact a \emph{K\"ahler--Einstein metric}, satisfying $\Ric\omega=-\lambda\omega$.
%
%
\subsection{The YTD conjecture}
An early obstruction for the existence of a cscK metric in $c_1(L)$ was found by Matsushima and Lichnerowicz: the identity component $\Aut^0(X,L)$ of the automorphism group of $(X,L)$ must be \emph{reductive}. This implies that the blowup of $\P^2$ at a point does not admit a cscK metric (for any $L$). This reductivity criterion is algebraic in nature, and an imprecise version of the Yau--Tian--Donaldson (YTD) conjecture states that the existence of a cscK metric in $c_1(L)$ is governed by a purely algebraic condition.
%
%
%
%
\section{Algebraic and non-Archimedean geometry}\label{sec:nageom}
As before, $X$ is a smooth complex algebraic variety, and $L$ a line bundle on $X$. The analytic geometry in~\S\ref{sec:archgeom} uses the standard absolute value on the ground field $\C$. In~\cite{BerkBook}, Berkovich developed a parallel analytic theory using the (non-Archimedean) \emph{trivial} absolute value on $\C$ (or any  field). We review some of its features here, following~\cite{trivval}.
%
%
\subsection{Divisorial valuations}\label{sec:divval}
First assume $X$ is projective. For any prime divisor $E$ over $X$ (\ie prime divisor on $X'$, for some smooth variety $X'$ with a projective birational morphism $X'\to X$), we can consider the order of vanishing $\ord_E\colon\C(X)^\times\to\Z$ along $E$. We define the set $X_\div$ of \emph{divisorial valuations} on $X$ as the set of valuations $v:=c\ord_E\colon\C(X)^\times\to\Q
$, where $E$ is a divisor over $X$, and $c\in\Q_{\ge0}$. ({\seb By convention, when $c=0$, $c\ord_E$ is equal to} the trivial valuation on $\C(X)$.)
We equip $X_\div$ with the weakest topology for which $v\mapsto v(f)$ is continuous for all $f\in\C(X)^\times$. 

Divisorial valuations admit a different geometric realization; see~\cite[\S4]{BHJ1} for details on what follows. By a \emph{test configuration} for $X$ we mean a normal scheme $\cX$ with a projective morphism $\cX\to\A^1$, a $\G_m$-action on $\cX$ lifting the standard action on $\A^1$, and a $\G_m$-equivariant isomorphism $\cX|_{\A^1\setminus\{0\}}\simto X\times(\A^1\setminus\{0\})$. To each irreducible component $E$ of the central fiber $\cX_0$, we define a valuation $v_E\colon\C(X)^\times\to\Q$ as the restriction of $\ord_E(\cX_0)^{-1}\ord_E$ from $\C(\cX)\simeq\C(X)(t)$ to $\C(X)$. Then $v_E\in X_\div$ is a divisorial valuation, and all divisorial valuations on $X$ appear in this way. 
%
%
\subsection{Berkovich analytification}
We continue to assume that $X$ is projective. The space $X_\div$ above admits a natural compactification, the \emph{Berkovich analytification} $X_\na$, which is compact (and Hausdorff). Just like $X_\div$, it is of algebro-geometric nature in the sense that it only uses the structure of $\C$ as an algebraically closed field.

As a set, $X_\na$ is the disjoint union, over all scheme-theoretic points $\xi\in X$, of the set of valuations $\kappa(\xi)^\times\to\R$ that are trivial on $\C$. Sending a point in $X_\na$ to the corresponding point $\xi$  defines a continuous \emph{support map} $X_\na\to X$.  There is also a natural \emph{center map} $X_\na\to X$, whose restriction to $X_\div$ sends a divisorial valuation $c\ord_E$ to the image of $E$ on $X$ when $c\in\Q_{>0}$, and to the generic point of $X$ when $c=0$. The center map has the curious property of being \emph{anticontinuous}: the preimage of a (Zariski) closed subset is open.

%
%
\subsection{The non-projective case}
The analytification $X_\na$ is defined even when $X$ is not projective, and $X\mapsto X_\na$ is functorial. If $X$ is quasi-projective, and $X\hto\bar X$ is an open immersion with $\bar X$ projective, then $X_\na$ is the open subset of $\bar X_\na$ consisting of all points whose support lies in $X$. It is a locally compact Hausdorff space. The set of points whose \emph{center} lies in $X$ is a compact subset $X_\beth\subset X_\na$. We have $X_\beth=X_\na$ iff $X$ is proper.

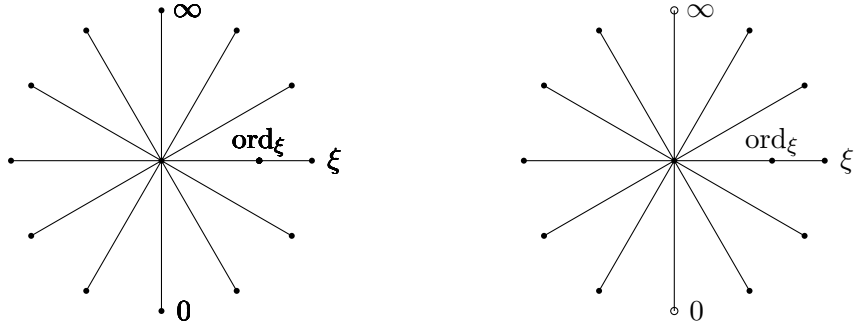
\begin{figure}
 \begin{tikzpicture}
    \tikzmath{\scpar = 0.5;}
    \tikzmath{\dang = 60;}
    \tikzmath{\x1 = 0; \y1 =0;}
    \tikzmath{\scale=2; \ang=0; }
    \tikzmath{\dotsize=0.8; }
    \tikzmath{\nobranches=12; }

    \node at (0,0) [circle,fill,inner sep=\dotsize]{};
    \foreach \i in {1,...,\nobranches} {%
      \tikzmath{\ang=360*\i/\nobranches;}
      \tikzmath{\x1 = \scale*cos(\ang); \y1 = \scale*sin(\ang);}
      \draw (0,0)--(\x1,\y1);
      \node at (\x1,\y1) [circle,fill,inner sep=\dotsize]{};
    \node at (1.3,0) [circle,fill,inner sep=\dotsize]{};
    \draw (0.35,2) node {$\infty$};
    \draw (0.3,-2) node {$0$};
    \draw (1.3,0.25) node {$\ord_\xi$};
   \draw (2.3,0) node {$\xi$};
    }
 \end{tikzpicture}
 \hspace*{20mm}
 \begin{tikzpicture}
    \tikzmath{\scpar = 0.5;}
    \tikzmath{\dang = 60;}
    \tikzmath{\x1 = 0; \y1 =0;}
    \tikzmath{\scale=2; \ang=0; }
    \tikzmath{\dotsize=0.8; }
    \tikzmath{\nobranches=12; }
    \tikzmath{\skipiA = int(\nobranches/4);}
    \tikzmath{\skipiB = int(3*\nobranches/4);}
    \tikzmath{\dotsizeopen = 1.2*\dotsize;}
    \node at (0,0) [circle,fill,inner sep=\dotsize]{};
    \foreach \i in {1,...,\nobranches} {%
      \tikzmath{\ang=360*\i/\nobranches;}
      \tikzmath{\x1 = \scale*cos(\ang); \y1 = \scale*sin(\ang);}
      \draw (0,0)--(\x1,\y1);
    \ifnum\i=\skipiA\relax
    \node at (\x1,\y1) [circle,draw,inner sep=\dotsizeopen]{};
    \else\ifnum\i=\skipiB\relax
    \node at (\x1,\y1) [circle,draw,inner sep=\dotsizeopen]{};
    \else
    \node at (\x1,\y1) [circle,fill,inner sep=\dotsize]{};
    \fi\fi
    }
    \node at (1.3,0) [circle,fill,inner sep=\dotsize]{};
    \draw (0.35,2) node {$\infty$};
    \draw (0.3,-2) node {$0$};
    \draw (1.3,0.25) node {$\ord_\xi$};
   \draw (2.3,0) node {$\xi$};
  \end{tikzpicture}
    \caption{The left figure shows the Berkovich space $\P^1_\na=\P^1_\beth$. The branch point is the trivial valuation on $\C(\P^1)$, and there is one branch for each closed point $\xi\in\P^1$. The interior of each branch is parametrized by $c\ord_\xi$, $0<c<\infty$, and the endpoint is the point $\xi$ itself (\ie the trivial valuation on $\kappa(\xi)\simeq\C)$. Any open neighborhood of the branch point contains all but finitely many full branches. The right figure shows $\G_{m,\na}=\P^1_\na\setminus\{0,\infty\}$. The space $\G_{m,\beth}$ is the compact subset of $\G_{m,\na}$ obtained by removing the open branches leading to $0$ and $\infty$. Finally, $\P^1_\div=\G_{m,\div}$ is the dense, totally disconnected subset of $\P^1_\na$ consisting of the branch point and all rational points $c\ord_\xi$, $c\in\Q_{>0}$.}
    \label{fig:Berkcurves}
\end{figure}

See Figure~\ref{fig:Berkcurves} for an illustration of $X_\na$ and $X_\beth$ for $X=\P^1$ and $X=\G_m$.

For any smooth curve $X$, the structure of $X_\na$ and $X_\beth$ is similar. In higher dimensions, the situation is more complicated, but some subsets can be understood. For example, the \emph{valuative tree}~\cite{valtree} and \emph{valuative tree at infinity}~\cite{eigval} are compact subsets of $\A^2_\an$. We refer to~\cite{dynberko} for a discussion, and to~\cite{hiro,izumi} for higher dimensions. These spaces can be used to study polynomial dynamics~\cite{eigval,dyncomp,DF21} as well as the singularities of plurisubharmonic functions~\cite{valpsh,valmul,hiro}.
%
%
\subsection{Other ground fields}\label{sec:otherfields}
The assignment $X\mapsto X_\na$ is a special case of an analytification functor that assigns to any variety $X$ over a complete valued field $k$ a topological space $X_\an$ (together with a structure sheaf, which we ignore here). In the case of $\C$ equipped with the usual (Archimedean) absolute value, $X_\an$ is the complex analytic space $X$. When instead $\C$ is equipped with the trivial absolute value, $X_\an$ becomes the space $X_\na$ above. 

Analytifications over more general non-Archimedean fields naturally show up as follows. Consider a map $f\colon Y\to X$ of complex algebraic varieties, with induced continuous map $f_\na\colon Y_\na\to X_\na$. A point $v\in X_\na$ can be viewed as a point $\xi\in X$ together with a valuation on the residue field $\kappa(\xi)$. Let $\cH(v)$ be the completion of $\kappa(\xi)$. Then the fiber $f_\na^{-1}(v)$ is naturally the analytification of the base change of $f^{-1}(\xi)$ with respect to $\kappa(\xi)\hto\cH(v)$.
%
%
\subsection{Hybrid spaces}\label{sec:hyb}
\begin{figure}
\includegraphics[width=0.2\textwidth]{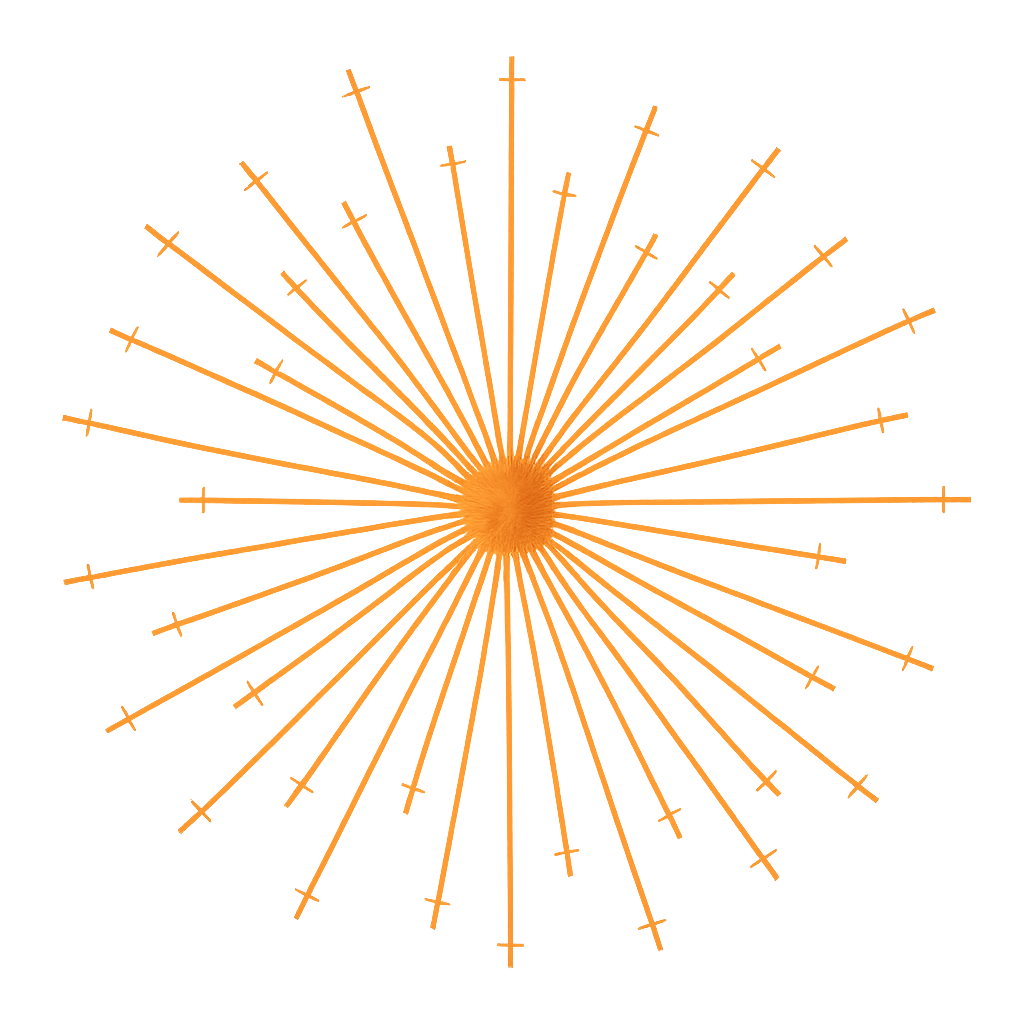}
\hspace*{10mm}
\includegraphics[width=0.2\textwidth]{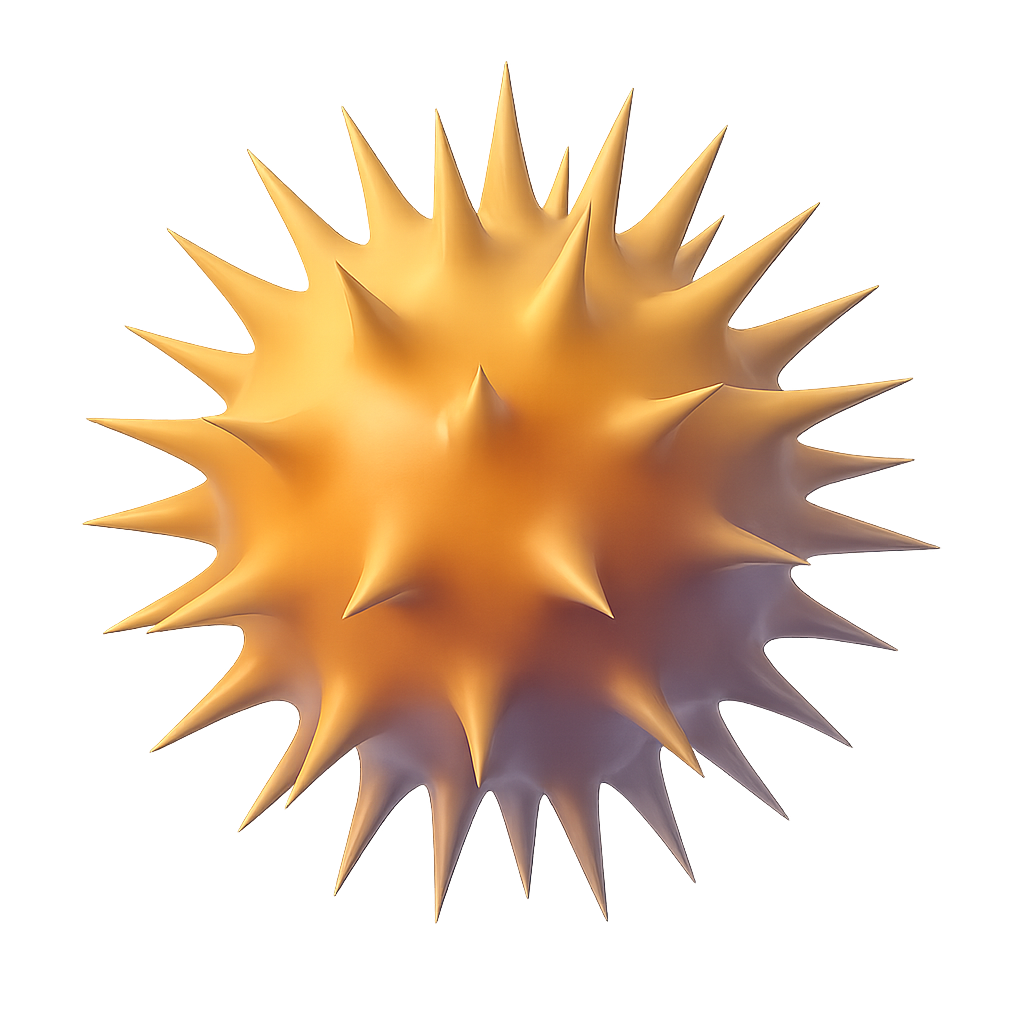}
\hspace*{10mm}
\includegraphics[width=0.2\textwidth]{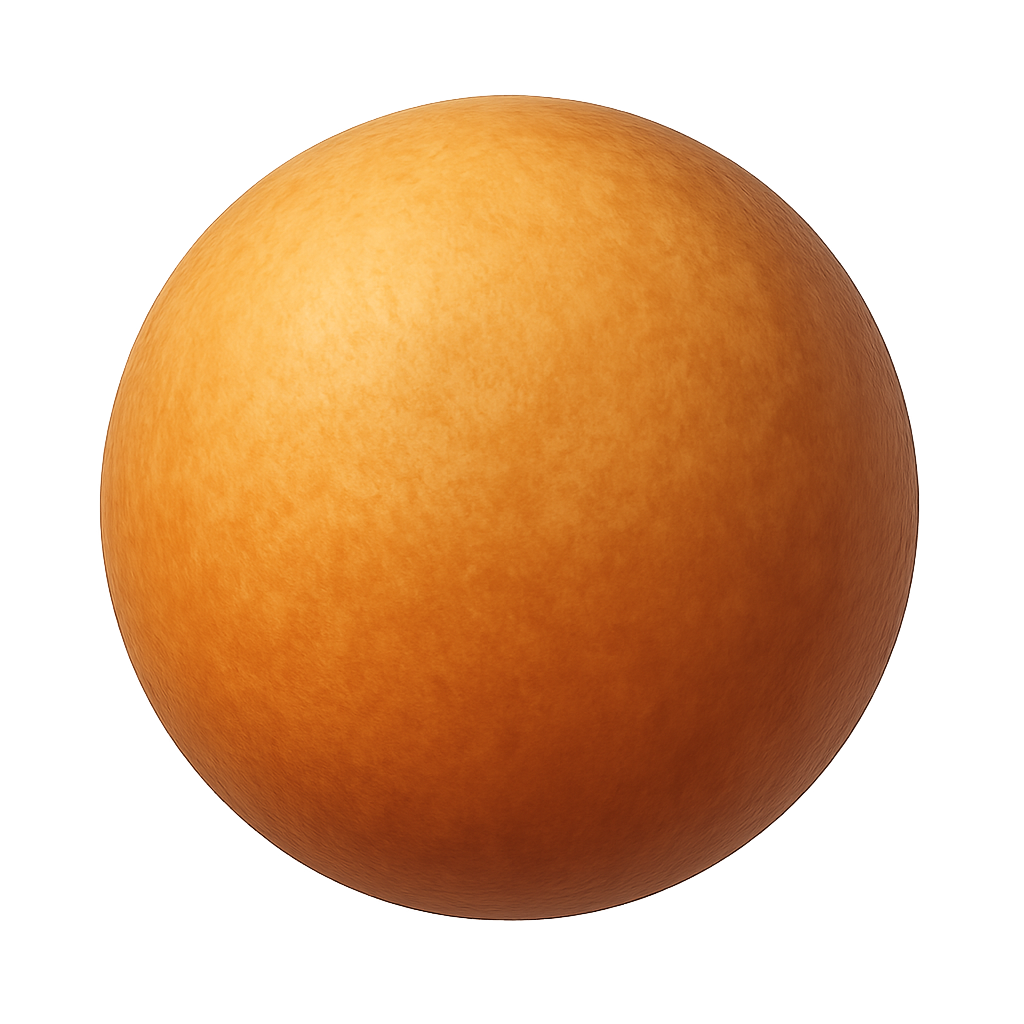}
\caption{The hybrid projective line $\P^1_\hyb$, fibering over the interval $[0,1]$. The round sphere on the right is the fiber over $t=1$ and equals the Riemann sphere. The spiky sphere in the middle illustrates the fiber over $t\in(0,1)$, which is homeomorphic to the Riemann sphere. The left-hand picture shows the fiber over $t=0$, the Berkovich projective line (see Figure~\ref{fig:Berkcurves}).}\label{fig:P1hyb}
\end{figure}

One can also consider analytifications over more general Banach rings. For example, using the field of complex numbers endowed with the maximum of the standard and trivial norms gives rise to \emph{hybrid spaces}, allowing us to degenerate complex manifolds to Berkovich spaces,~\cite{Berk09,konsoib}. Specifically, a smooth complex algebraic variety $X$ gives rise to a locally compact topological space $X_\hyb$ fibering over the interval $[0,1]$ (the hybrid point). The fiber over $t\in(0,1]$ is the analytification of $X$ with respect to the absolute value $|\cdot|_\infty^t$ on $\C$, and is homeomorphic to the complex manifold $X$. The fiber over $t=0$ equals the Berkovich space $X_\na$. See Figure~\ref{fig:P1hyb} for an illustration of the hybrid projective line. While we do not explicitly use them in the present work, hybrid spaces can be used to understand and formulate some of the constructions that we do employ.

%
%
\subsection{Metrics on line bundles}
Let $X$ be an algebraic variety and $L$ a line bundle on $X$. Let $X_\na$ and $L_\na$ be the analytifications of $X$ and (the total space of) $L$, respectively. Denote by $L_\na^\times$ the analytification of the total space of $L$ with the zero section removed.

We have a continuous surjective map $L_\na^\times\to X_\na$. The fiber $L_\na^\times(v)$ over a point $v\in X_\na$ is isomorphic to the analytification of the multiplicative group $\G_{m,\cH(v)}$ over the non-Archimedean field $\cH(v)$; see~\S\ref{sec:otherfields}.

By a \emph{metric} (resp.\ \emph{singular metric}) on $L_\na$ we mean a function $\phi\colon L_\na^\times\to\R$ (resp.\ $\phi\colon L_\na^\times\to\R\cup\{\pm\infty\}$) with the same log-homogeneity property along fibers as in the complex analytic case; this property implies that $\phi|_{L_\na^\times}$ is determined by its value at any point. 

The same operations as in~\S\ref{sec:ametrics} are available in the current non-Archimedean setting. Any section $s\in\Hnot(X,L)$ defines a continuous map $s\colon X_\na\to L_\na$, and a singular metric $\phi:=\log|s|$, uniquely defined by $\phi\circ s(v)=0$ when $s(v)\ne0$, and $\phi\equiv-\infty$ on $L_\na^\times(v)$ otherwise. The metric $\log|s|$ is continuous where $s\ne 0$.

It is sometimes convenient to say "non-Archimedean metric on $L$" instead of "metric on $L_\na$", and we shall do so occasionally.
%
%
\subsection{Fubini--Study metrics}\label{sec:naFS}
Assume $X$ is projective and $L$ is ample. An analog of the space $\cH$ in~\S\ref{sec:ametrics} is the space $\cH_\na=\cH_\na(L)$ of \emph{non-Archimedean Fubini--Study metrics} on $L$, that is, the space of metrics on $L_\na$ of the form 
\begin{equation}\label{equ:FSna1}
  \phi=\frac1m\max_j(\log|s_j|+\la_j),
\end{equation}
where $m\ge 1$, $\{s_j\}_j$ is a finite set of global sections of $mL$ without common zeros, and $\la_j\in\Z$. Any Fubini--Study metric is continuous, so the difference of any two Fubini--Study metrics on $L$ is a continuous function on $X_\na$. In fact, the set of such differences is dense in $C^0(X_\na)$. 

The Berkovich space $X_\na$ is not a differentiable manifold, even when $X$ is smooth. Chambert--Loir and Ducros~\cite{CLD25} nevertheless defined a notion of a smooth metric, and developed a theory of forms and currents on Berkovich spaces. We will not directly rely on their work, but the global pluripotential theory outlined in~\S\ref{sec:gppt} is compatible with theirs. A non-Archimedean Fubini--Study metric is not smooth, but can be approximated by smooth ones. 

%
%
\subsection{The trivial metric}
When $X$ is projective (or proper), any line bundle $L$ on $X$ admits a canonical non-Archimedean metric, the \emph{trivial metric} $\phi_\triv$. This means that Calabi's question has a trivial answer in the trivially valued case!

The trivial metric is continuous, and defined as follows: pick any $v\in X_\na$, with center $\eta\in X$, and let $s$ be a trivializing section of $L$ in an open neighborhood of $\eta$. We then define $\phi_\triv$ on $L^\times_\na(v)$ by $\phi_\triv(s(v))=0$. This definition does not depend on the choice of trivializing section.

When $L$ is ample, the trivial metric is a Fubini--Study metric of the following type. Pick $m\ge1$ and a finite set $(s_j)$ of nonzero global sections of $mL$ without common zeros. Then $\phi_\triv=\frac1m\max_j\log|s_j|$. 

The trivial metric allows us to identify metrics on $L_\na$ with functions on $X_\na$. For example, if $s\in H^0(X,L)$, then $\log|s|-\phi_\triv$ is a singular metric on $\cO_{X,\na}$, \ie a function on $X_\na$ whose value at $v\in X_\na$ we write as $-v(s)$, where $v(s)\in [0,+\infty]$. To a Fubini--Study metric $\phi$ as in~\eqref{equ:FSna1} is associated a \emph{Fubini--Study potential} $\f=\phi-\phi_\triv$, a continuous function on $X_\na$ given by 
\[
\f(v)=\frac1m\max_j(\la_j-v(s_j)).
\]
The potential associated to the trivial metric is the zero function.
%
%
\subsection{Test configurations}\label{sec:tc}
Assume $X$ is projective, and consider a test configuration $(\cX,\cL)$ for $(X,L)$. This means that $\cX$ is a test configuration for $X$ as in~\S\ref{sec:divval}, and that $\cL$ is a $\Q$-line bundle such that $\cL|_{\cX\setminus\cX_0}$ is $\G_m$-equivariantly isomorphic to $p_1^*L$ via $\cX\setminus\cX_0\simeq X\times(\A^1\setminus\{0\})$. Assuming that $\cX$ is normal, $(\cX,\cL)$ defines a function $\f=\f_{\cX,\cL}$ on $X_\div$ as follows. Given any divisorial valuation $v\in X_\div$, there exists a normal (even smooth) test configuration $\cX'$ for $X$, such that $v$ is defined by an irreducible component $E\subset\cX'_0$ and the canonical $\G_m$-equivariant birational maps $\rho\colon\cX'\dashrightarrow\cX$ and $\mu\colon\cX'\dashrightarrow X\times\A^1$ are morphisms. Then we can write $\rho^*\cL-\mu^*\mathrm{pr_1}^*L=\cO_{X'}(D)$ for a $\Q$-Cartier $\Q$-divisor supported on $\cX'_0$, and we set $\f(v)=\ord_E(D)/\ord_E(\cX'_0)$.

The function $\f=\f_{\cX,\cL}$ extends to a continuous function on $X_\na$. Moreover, $\f_\na$ is a Fubini--Study potential when $(\cX,\cL)$ is ample, in the sense that $\cL$ is relatively ample. In fact $(\cX,\cL)\mapsto\f_{\cX,\cL}$ gives a bijection between the set of normal ample test configurations and $\cH_\na$. Here the trivial test configuration $(X,L)\times\A^1$ corresponds to the zero function.

It is sometimes convenient to think of a test configuration $(\cX,\cL)$ as defined over $\P^1$ rather than $\A^1$: this can be done by gluing together $(\cX,\cL)$ to the trivial test configuration.
%
%
\subsection{Log discrepancy and Temkin's metric}\label{sec:Temkin}
Assume finally that $X$ is smooth and projective. Following the Minimal Model Program in birational geometry, we define the \emph{log discrepancy function} $A_X\colon X_\div\to\Q_{\ge 0}$ as the unique function that is homogeneous with respect to the scaling action of $\Q_{>0}$ and such that 
$$
A_X(\ord_E)=1+\ord_E(K_{Y/X})
$$
for any prime divisor $E$ on a smooth birational model $Y\to X$. By~\cite{hiro,JM12}, the log discrepancy function is lsc on $X_\div$, and hence admits a (unique) maximal lower-semicontinuous extension $A_X\colon X_\na\to [0,+\infty]$. 

Using the trivial metric on the canonical bundle $K_X$, the log discrepancy function $A_X$ can be viewed as a singular non-Archimedean metric on $K_X$, which turns out to be compatible with the general theory of Temkin~\cite{Tem16}. In the complex analytic case, this metric would correspond to a positive measure on $X$. In the present non-Archimedean context, this is not the case anymore, but the function $e^{-2A_X}$ on $X_\na$ can nevertheless be interpreted, for some purposes, as the density of this would-be measure on $X_\na$.
%
%
%
%
\section{Global pluripotential theory and cscK metrics}\label{sec:gppt}
As above, fix a polarized smooth projective complex variety $(X,L)$. The study of plurisubharmonic (psh) metrics on $L$, viewed as a holomorphic line bundle, is well developed and forms part of global (or geometric) pluripotential theory; see for example~\cite{GZ,BBGZ,BBEGZ}. In~\cite{trivval,nakstab2,nakstab1} we showed that many of the basic features have non-Archimedean analogs, the ground field $\C$ being equipped with the trivial valuation. Here we summarize the main results for the two cases. Note that the presentation does not always follow the order in which the results are actually proved.
%
%
\subsection{Plurisubharmonic metrics}
As explained above, any global section $s\in\Hnot(X,mL)\setminus\{0\}$, $m\in\Z_{>0}$, defines a singular Hermitian metric $m^{-1}\log|s|$ on $L$. The set $\PSH=\PSH(L)$ of psh metrics on $L$ can be characterized as the smallest set of singular metrics on $L$, not identically $-\infty$, that contains all metrics induced by sections, and is further invariant under addition of constants, finite maxima, and decreasing limits. The latter means that if $(\phi_j)$ is a decreasing sequence in $\PSH$, then $\phi:=\lim_j\phi_j$ either lies in $\PSH$, or satisfies $\phi\equiv-\infty$.

The choice of a continuous reference metric $\phi_\refe$ on $L$ allows us to view the elements of $\PSH$ as functions, or \emph{potentials} on $X$, via $\phi\mapsto\phi-\phi_\refe$. These functions are upper semicontinuous and integrable, and the \emph{weak topology} on $\PSH$ can be defined as the induced $L^1$-topology. 

Similarly, any $s\in\Hnot(X,mL)\setminus\{0\}$ also defines a singular non-Archimedean metric $m^{-1}\log|s|$ on $L_\na$, and the set $\PSH_\na=\PSH_\na(L)$ of \emph{non-Archimedean psh metrics} on $L$ is characterized in exactly the same way as above. The trivial metric on $L_\na$ is a canonical choice of reference metric, and allows us to view the elements of $\PSH_\na$ as usc functions $\f\colon X_\na\to\R\cup\{-\infty\}$. Their restrictions to $X_\div\subset X_\na$ are in fact finite-valued, and the weak topology on $\PSH_\na$ is the one given by pointwise convergence on $X_\div$.

A general metric in $\PSH$ (resp.~$\PSH_\na$) is usc, and possibly unbounded; we denote by $\CPSH$ (resp.\ $\CPSH_\na$) the subset of continuous (bounded) psh metrics.

Any smooth positive metric $\phi\in\cH$ is psh, and every element of $\PSH$ (resp.\ $\CPSH$) is a decreasing (resp.\ uniform) limit of elements in $\cH$. On the \nA side, the same results hold for $\cH_\na$, the set of non-Archimedean Fubini--Study metrics. 

The spaces $\PSH$ and $\PSH_\na$ enjoy further deep properties, notably the \emph{continuity of envelopes}: if $\phi$ is a continuous metric on $L$ (resp.\ $L_\na$), then there is a (unique) largest psh metric $\pp(\phi)$ such that $\pp(\phi)\le\phi$, and $\pp(\phi)$ lies in $\CPSH$ (resp.\ $\CPSH_\na$). 

We will not explore this here, but it is also possible to partially develop global pluripotential theory on hybrid spaces (see~\S\ref{sec:hyb}), combining the complex analytic and trivially valued worlds, see~\cite{PS23,Reb23,Hybrid}.
%
%
\subsection{Monge--Amp\`ere operator and energy functionals}
The \emph{complex Monge--Amp\`ere operator} takes a metric $\phi\in\cH$ to the smooth probability measure 
$$
\MA(\phi):=V^{-1}(\ddc\phi)^n,\quad V=(L^n),
$$
\ie the volume form of the Kähler metric $\ddc\phi$, normalized to mass $1$. This operator admits an \emph{Euler--Lagrange functional}, \ie a functional $\ee\colon\cH\to\R$ whose directional derivatives at $\phi\in\cH$ are obtained by integration against $\MA(\phi)$. It is known as the \emph{Monge--Amp\`ere energy}, and is explicitly given by 
\[
\ee(\phi)=\frac1{n+1}\sum_{j=0}^n V^{-1}\int_X(\phi-\phi_\refe)(\ddc\phi)^j\wedge(\ddc\phi_\refe)^{n-j},
\]
after normalization by $\ee(\phi_\refe)=0$, where $\phi_\refe\in\cH$ is a fixed reference metric.
A related functional is the \emph{Ricci energy} $\rr\colon\cH\to\R$, given by the similar looking expression
\[
\rr(\phi):=\sum_{j=0}^{n-1}V^{-1}\int_X(\phi-\phi_\refe)\,\theta\wedge(\ddc\phi)^j\wedge(\ddc\phi_\refe)^{n-j-1}, 
\]
where $\theta:=\ddc\log(\ddc\phi_\refe)^n$ is the curvature form of the induced reference metric on the canonical bundle $K_X$, \ie minus the Ricci curvature of the reference Kähler metric $\ddc\phi_\refe$. The Ricci energy again arises as the Euler--Lagrange functional of the measure-valued operator $\cH\ni\phi\mapsto n V^{-1}\theta\wedge(\ddc\phi)^{n-1}$. 

To make the analogous definitions in the non-Archimedean case, we use the fact that each $\f\in\cH_\na$ uniquely corresponds to a normal ample test configuration $(\cX,\cL)\to\P^1$; the \emph{non-Archimedean Monge--Amp\`ere energy} of $\f$ can then be written as the normalized top self-intersection number 
\[
\ee_\na(\f)=\frac{(\cL^{n+1})}{(n+1)V}. 
\]
The \emph{non-Archimedean Ricci energy} $\rr_\na(\phi)$ can similarly be described in terms of intersection numbers, as can the \emph{non-Archimedean Monge--Amp\`ere operator}

$$
\MA_\na\colon\cH_\na\to\cM_\na, 
$$
which computes the directional derivatives of $\ee_\na\colon\cH_\na\to\R$. For each $\f\in\cH_\na$, $\MA_\na(\f)$ is a \emph{divisorial measure}, \ie a probability measure on $X_\na$ with support a finite set of divisorial valuations. Concretely, if $\f\in\cH$ is defined by a normal ample test configuration $(\cX,\cL)$, with central fiber $\cX_0=\sum_ib_iE_i$, then 
\[
\MA_\na(\f)=\sum_ib_i(\cL|_{E_i}^n)v_i,
\]
where $v_i\in X_\div$ is the divisorial valuation associated to $E_i$, see~\ref{sec:divval}.

%
%
\subsection{Metrics of finite energy}\label{sec:E1}
The Monge--Amp\`ere energy $\ee\colon\cH\to\R$ is usc and increasing, so its minimal usc extension $\ee\colon\PSH\to\R\cup\{-\infty\}$ is given by 
\[
\ee(\phi)=\inf\{\ee(\psi)\mid \cH\ni\psi\ge\phi\}.
\]
We define the space $\cE^1\subset\PSH$ of \emph{metrics of finite energy} as the set
\[
\cE^1:=\{\phi\in\PSH\mid\ee(\phi)>-\infty\}, 
\]
and equip it with the \emph{strong topology}, the coarsest refinement of the weak topology in which $\ee\colon\cE^1\to\R$ becomes continuous.

In the non-Archimedean case, we define $\ee_\na\colon\PSH_\na\to\R\cup\{-\infty\}$, the subset $\cE^1_\na\subset\PSH_\na$, and the strong topology in exactly the same way.

Any bounded psh metric is of finite energy, so we have strongly dense inclusions $$
\cH\subset\CPSH\subset\cE^1,\quad\cH_\na\subset\CPSH_\na\subset\cE^1_\na.
$$
By definition, $\ee\colon\cE^1\to\R$ is strongly continuous. The Ricci energy and Monge--Amp\`ere operator also admit (unique) strongly continuous extensions
$$
\rr\colon\cE^1\to\R,\quad\MA\colon\cE^1\to\cM,
$$
where $\cM$ is the space of probability measures on $X$, equipped with the weak topology. The analogous statements hold on the non-Archimedean side.

We note that $\cE^1_\na$ contains many objects familiar to algebraic geometers, for example quite general filtrations of the section ring of $L$.
%
%
\subsection{The Calabi--Yau theorem}
Define the (Monge--Amp\`ere) \emph{energy} of a probability measure $\mu\in\cM$ by 
\[
\ee^\vee(\mu):=\sup_{\phi\in\cH}\left\{\ee(\phi)-\int(\phi-\phi_\refe)\,\mu\right\}\in [0,+\infty]. 
\]
We have $\ee^\vee(\mu)\ge0$, with equality iff $\mu=\mu_\refe:=\MA(\phi_\refe)$. The energy functional $\ee^\vee\colon\cM\to[0,+\infty]$ is convex and lower semicontinuous (lsc) in the weak topology, and we equip the set 
\[
\cM^1:=\{\mu\in\cM
\mid\ee^\vee(\mu)<\infty\}
\]
of measures of finite energy with the strong topology, the coarsest refinement of the weak topology in which $\ee^\vee\colon\cM^1\to\R$ is continuous. Both the set $\cM^1$ and its strong topology turn out to be independent of the choice of $L$.

A finite energy version of the Calabi--Yau theorem now says that
\[
\MA\colon\cE^1/\R\to\cM^1
\]
is a homeomorphism in the strong topology. Given a measure $\mu\in\cE^1$, the solution $\phi\in\cE^1$ to $\MA(\phi)=\mu$, unique up to a constant, is precisely the maximizer of $\ee(\phi)-\int(\phi-\phi_\refe)\mu$. 

All these results have exact analogs on the non-Archimedean side: we define $\ee_\na^\vee$, $\cM^1_\na$ in the same way, and obtain a homeomorphism 
$$
\MA_\na\colon\cE^1_\na/\R\simto\cM^1_\na
$$
in the strong topology. However, things do diverge somewhat when it comes to regularity of solutions. The fundamental theorem of Yau~\cite{Yau78} states that $\phi\in\cE^1$ lies in $\cH$ iff $\MA(\phi)$ is a (smooth, positive) volume form. On the non-Archimedean side, we have seen that $\MA_\na(\phi)$ is a divisorial measure for each $\phi\in\cH_\na$; conversely, each $\phi\in\cE^1$ with $\MA_\na(\phi)$ a divisorial measure lies in $\CPSH_\na$, but a precise characterization of when $\phi$ actually lies in $\cH_\na$ seems currently out of reach. 
%
%
\subsection{Entropy}
The \emph{entropy} of a probability measure $\mu\in\cM$ on the complex manifold $X$ is defined as its relative  entropy\footnote{The factor one-half is included to ensure compatibility with the non-Archimedean case below.} with respect to the reference measure $\mu_\refe=\MA(\phi_\refe)$, \ie 
$$
\Ent(\mu):=\frac 12\int_X\log\left(\frac{\mu}{\mu_\refe}\right)\mu\in [0,+\infty]
$$
if $\mu$ is absolutely continuous with respect to $\mu_\refe$, and $\Ent(\mu)=+\infty$ otherwise. This defines a convex, lsc functional 
$$
\Ent\colon\cM\to[0,+\infty],
$$
which satisfies $\ee^\vee(\mu)\le C\left(\Ent(\mu)+1\right)$ for a uniform constant $C>0$. In particular, every measure of finite entropy has finite energy. A crucial property is that sets of bounded entropy are further \emph{compact} for the strong topology of $\cM^1$.

\medskip

Consider next the non-Archimedean case, and recall that the function $e^{-2A_X}$ on the Berkovich space $X_\na$ can be thought of as the density of a would-be measure on $X_\na$ corresponding to the singular non-Archimedean metric on $K_X$ induced by the log discrepancy function $A_X\colon X_\na\to [0,+\infty]$ (see~\S\ref{sec:Temkin}). With this in mind, we define the \emph{non-Archimedean entropy} of a measure $\mu\in\cM_\na$ as 
$$
\Ent_\na(\mu):=\int_{X_\na} A_X\,\mu\in [0,+\infty]. 
$$
This yields again an lsc functional $\Ent_\na\colon\cM_\na\to [0,+\infty]$ that satisfies\footnote{Here the additive constant disappears by homogeneity with respect to the scaling action of $\R_{>0}$ on $X_\na$.} 
$$
\ee_\na^\vee(\mu)\le C\Ent_\na(\mu)
$$
for a uniform constant $C>0$. However, a key difference is that sets of bounded entropy are no longer strongly compact in $\cM^1_\na$, exhibiting one of the relatively few instances where the two worlds diverge.

\subsection{CscK metrics and the Mabuchi functional}\label{sec:cscK}
The scalar curvature of the K\"ahler metric $\ddc\phi$ defined by $\phi\in\cH$ is given by 
\[
S(\phi)=-n\frac{\ddc\log(\ddc\phi)^n\wedge(\ddc\phi)^{n-1}}{(\ddc\phi)^n}\in C^\infty(X),
\] 
and we say that $\ddc\phi$ is a cscK metric if $S(\phi)$ is constant, necessarily equal to 
\[
\bar S:=-n(K_X\cdot L^{n-1})/(L^n).
\]
Thus cscK metrics correspond to zeros of the operator 
$$
\phi\mapsto(\bar S-S(\phi))\MA(\phi)=\bar S\MA(\phi)+n V^{-1}\ddc\log(\ddc\phi)^n\wedge(\ddc\phi)^{n-1},
$$
which also admits an Euler--Lagrange functional
$$
\mm\colon\cH\to\R,
$$
known as the \emph{Mabuchi K-energy} and explicitly given by the Chen--Tian formula
\begin{equation}\label{equ:CT}
\mm(\phi)=\bar S\ee(\phi)+\rr(\phi)+\hh(\phi),\quad\text{where}\quad\hh(\phi):=\Ent(\MA(\phi)).
\end{equation}
This formula provides a natural lsc extension $\mm\colon\cE^1\to\R\cup\{+\infty\}$, which can be characterized as the maximal lsc extension from $\cH$. In view of Yau's theorem, this boils down to the fact that any $\mu\in\cM^1$ can be written as the strong limit of volume forms $\mu_j$ such that $\Ent(\mu_j)\to\Ent(\mu)$, which can be proved by a standard regularization argument. 

\medskip
The \emph{non-Archimedean Mabuchi functional} 
$$
\mm_\na\colon\cE^1_\na\to\R\cup\{+\infty\}
$$
is similarly defined by 
\begin{equation}\label{equ:CTna}
\mm_\na=\bar{S}\ee_\na+\hh_\na+\rr_\na,\quad\text{where}\quad\hh_\na(\phi):=\Ent_\na(\MA_\na(\phi)).
\end{equation}
While it is also lsc, a key difference with the complex case is that it is not known to be the \emph{maximal} lsc extension from $\cH_\na$. The \emph{entropy regularization conjecture} posits that this holds, and is equivalent to the existence, for every $\mu\in\cM^1_\na$ on $X_\na$, of a sequence $(\mu_j)$ in the image of $\MA_\na\colon\cH_\na\to\cM^1_\na$ such that $\mu_j\to\mu$ strongly in $\cM^1_\na$ and 
$$
\Ent_\na(\mu_j)=\int_{X_\na} A_X\,\mu_j\to\Ent_\na(\mu)=\int_{X_\na} A_X\,\mu.
$$
Using known results on non-Archimedean Monge--Amp\`ere equations, this can be reduced to the case where $\mu$ is a divisorial measure, but the major stumbling block is then the lack of characterization of $\MA_\na(\cH_\na)$ within divisorial measures, already mentioned above. 
%
%
%
%
\section{Finsler geometry of finite energy metrics}\label{sec:E1geom}
To find a cscK metric, the idea is to use a variational method, and study the Mabuchi functional $\mm$ on the space $\cH$ of smooth positive metrics on $L$: indeed, as explained above cscK metrics are exactly the curvature forms of critical points of $\mm$. To this end, we introduce a Finsler metric on $\cH$ whose completion is a geodesic metric space with respect to which the Mabuchi functional turns out to be convex.

%
%
\subsection{Finsler metrics on the space of Kähler potentials}
The space $\cH$ sits as an open convex subset of the affine space $C^\infty(L)$ of smooth Hermitian metrics on $L$, and its tangent space at $\phi\in\cH$ can thus be canonically identified with the vector space $C^\infty(X)$ of smooth functions. 

For any $p\in [1,\infty]$, the tangent bundle of $\cH$ is thus equipped with a \emph{tautological $L^p$-norm}, restricting to the $L^p(\MA(\phi))$-norm on each tangent space $T_\phi\cH\simeq C^\infty(X)$. When $p=2$, this defines a Riemannian metric on $\cH$, first considered by Mabuchi, who further proved that the Mabuchi energy $\mm\colon\cH\to\R$ is \emph{infinitesimally convex} in this metric, in the sense that its Hessian is everywhere semipositive. 

While the infinite dimensional Riemannian manifold $\cH$ does not admit geodesics in general~\cite{LV13}, X.X.Chen~\cite{Che00} was nevertheless able to construct weak geodesics, obtained as limits of smooth `almost' geodesics in $\cH$. This allowed him to show that the associated length pseudometric $\dd_2$ is an actual metric (\ie separates points), something that was later extended by Darvas to the length metric $\dd_p$ of each tautological $L^p$-norm, $p\in [1,\infty]$. 

In fact, it is the case $p=1$ that turned out to have most relevance to  pluripotential theory. First, the $\dd_1$ metric on $\cH$ can be characterized in terms of the Monge--Amp\`ere energy $\ee\colon\cH\to\R$, as it satisfies, for all $\phi,\p\in\cH$,
\begin{equation}\label{equ:d1}
\dd_1(\phi,\p)=\ee(\phi)+\ee(\p)-2\sup_{\cH\ni\tau\le\phi,\p}\ee(\tau).
\end{equation}
Answering a conjecture of Guedj, Darvas proved that the metric completion of $(\cH,\dd_1)$ can be identified with the space $\cE^1$ of psh metrics of finite energy with its strong topology, see~\S\ref{sec:E1}. For any $\phi,\psi\in\cE^1$, he further showed that the set of metrics in $\cE^1$ dominated by both $\phi$ and $\psi$ admits a maximal element $\pp(\phi,\psi)\in\cE^1$, and~\eqref{equ:d1} becomes
\begin{equation}\label{equ:d1alt}
\dd_1(\phi,\psi)=\ee(\phi)+\ee(\psi)-2\ee(\pp(\phi,\psi)). 
\end{equation}

As proved in~\cite{nakstab1}, this can be adapted to the non-Archimedean case by adopting the analog of~\eqref{equ:d1} as a \emph{definition} of $\dd_{1,\na}$ on $\cH_\na$. One then proves that $\dd_{1,\na}$ indeed defines a metric (the triangle inequality being the delicate part), and that the completion of $(\cH_\na,\dd_{1,\na})$ equals the space $\cE^1_\na$ of metrics of finite energy, the analog of~\eqref{equ:d1alt} also holding true in that case. 
%
%
\subsection{Busemann convexity and maximal geodesic rays}
Chen's weak geodesics mentioned above ultimately yield a distinguished class of \emph{psh} geodesics in the complete metric space $(\cE^1,\dd_1)$, with respect to which it is uniquely geodesic, and even Busemann convex. 

Just as in~\S\ref{sec:radial}, we can thus define $\cE^1_\ra$ as the set of directions of psh geodesic rays $\{\phi_t\}_{t\ge 0}$ in $(\cE^1,\dd_1)$, which becomes a complete metric space with respect to the radial metric
$$
\dd_{1,\ra}(\{\phi_t\},\{\psi_t\}):=\lim_{t\to\infty}t^{-1}\dd_1(\phi_t,\psi_t).
$$

In analogy with~\S\ref{sec:nanorms}, we can further define two maps 
$$
\cE^1_\ra\twoheadrightarrow\cE^1_\na,\quad\cE^1_\na\hto\cE^1_\ra,
$$
the difference being that the latter map is only right-inverse to the former in the present infinite-dimensional context. Let us briefly describe these two maps. 

First consider a psh geodesic ray $\{\phi_t\}\in\cE^1_\ra$, to which we then want to associate a function $\f\in\cE^1_\na$.  If $\om_\refe:=\ddc\phi_\refe$, then $\{\phi_t\}$ can be viewed as an $S^1$-invariant $\Omega$-psh function $\Phi$ on $X\times\DD^\times$ via $\Phi(x,\tau)=(\phi_{-\log|\tau|}-\phi_\refe)(x)$; here $\Om$ is the pullback of $\om$ to $X\times\DD^\times$.
By~\cite[Theorem~1]{Dar17}, there exists $A>0$ such that $\tilde\Phi:=\Phi+A\log|\tau|$ extends to an $\Om$-psh function on $X\times\DD$. If $\cX$ is a smooth test configuration for $X$ that dominates $X\times\P^1$ in the sense that the 
canonical birational map $\rho\colon\cX\dashrightarrow X\times\P^1$ is a morphism, then $\tilde\Phi$ also extends to a $\rho^\star\Om$-psh function on $\cX|_\DD$. Any irreducible component $E$ of the central fiber $\cX_0=\sum_Eb_EE$ induces a divisorial valuation $v_E\in X_\div$, see~\S\ref{sec:divval}, and we set 
\[
\f(v_E)=-b_E^{-1}\ord_E(\tilde\Phi)-A,
\] 
where $\ord_E(\tilde\Phi)$ is the (generic) Lelong number of $\tilde\Phi$ along $E$. Now, any divisorial valuation on $X$ is of the form $v_E$ for some irreducible component $E$ of some smooth test configuration $\cX$ for $X$ dominating $X\times\P^1$, and it is not hard to see that the construction above defines a function $\f\colon X_\div\to\R$ that depends only on $\{\phi_t\}\in\cE^1_\ra$, and not on any choices made. It was proved in~\cite{YTD} (using multiplier ideals) that $\f$ extends (uniquely) to a function $\f\in\PSH_\na$, and that $\f\in\cE^1_\na$; in fact $\ee_\na(\f)\ge\lim_{t\to\infty}\ee(\phi_t)$.

To describe the map $\cE^1_\na\hto\cE^1_\ra$, first consider a non-Archimedean Fubini--Study metric $\f\in\cH_\na$. This is associated with a normal, ample test configuration $(\cX,\cL)$ for $(X,L)$. Berman~\cite{Berm16} used an envelope construction on $\cX$ to produce a psh geodesic ray $\{\phi_t\}_{t\ge0}$ in $\cE^1$ emanating from $\phi_\refe$. It was proved in~\cite{YTD} that the resulting map $\cH_\na\to\cE^1_\ra$ is injective, and that 
\begin{equation}\label{equ:eeasym}
\lim_{t\to\infty} t^{-1}\ee(\phi_t)=\ee_\na(\f);
\end{equation}
moreover, $\f\in\cH_\na\subset\cE^1_\na$ is the image of $\{\phi_t\}\in\cE^1_\ra$ under the map $\cE^1_\ra\to\cE^1_\na$ above.

By approximation, we get the desired injection $\cE^1_\na\hto\cE^1_\ra$, which was subsequently proved to be an isometry by Reboulet~\cite{Reb23}. The geodesic rays in $\cE^1_\ra$ in the image of $\cE^1_\na$ are called \emph{maximal}. The same techniques for proving these results also show that various energy functionals have the desired asymptotics along maximal geodesic rays: in particular, 
\begin{equation}\label{equ:rrasym}
\lim_{t\to\infty} t^{-1}\rr(\phi_t)=\rr_\na(\f). 
\end{equation}

%
%
\subsection{Asymptotics of the Mabuchi functional}\label{sec:entasym}
A key result in~\cite{wytd} is now the following result on the asymptotics of the Mabuchi functional: for any maximal psh geodesic ray $\{\phi_t\}$ in $\cE^1$ with direction $\f\in\cE^1_\ra$, we have
\begin{equation}\label{equ:mabasym}
\lim_{t\to\infty} t^{-1}\mm(\phi_t)=
\begin{cases}
    \mm_\na(\f) &\text{if $\f\in\cE^1_\na$}\\
    +\infty &\text{otherwise}.
\end{cases}
\end{equation}
In view of~\eqref{equ:eeasym}, \eqref{equ:rrasym} and the Chen--Tian formula~\eqref{equ:CT}, \eqref{equ:CTna}, the difficulty here is to control the slope at infinity of the entropy functional $\hh(\phi_t)=\Ent(\MA(\phi_t))$. 

The crucial fact that $\hh(\phi_t)$ (and hence $\mm(\phi_t)$) has infinite slope when $\f\notin\cE^1_\na$ was proved by C.~Li~\cite{LiGeod}, who also showed 
$$
\lim_{t\to\infty} t^{-1}\hh(\phi_t)\ge\hh_\na(\f)
$$
when $\f\in\cE^1_\na$, building on~\cite{BHJ2}. Our proof of equality also leverages on~\cite{BHJ2} together with a key regularity result for psh envelopes in~\cite{BD12,DNT24}.

%
%
%
%
\section{The Yau--Tian--Donaldson conjecture}\label{sec:ytd}
We are now ready to sketch a proof of Theorem~A in the introduction.
%
%
\subsection{Stability notions}
The notions below are in terms of the non-Archimedean Mabuchi functional $\mm_\na\colon\cE^1_\na\to\R\cup\{+\infty\}$, but we emphasize that this means they are of algebro-geometric nature.
\begin{defi}\label{defi:hKstab}
    We say that $(X,L)$ is 
    \begin{itemize}
        \item[(i)] \emph{$\hK$-stable} if $\mm_\na(\f)\ge0$ for all $\f\in\cE^1_\na$, with equality iff $\f\in\R$ is constant;
        \item[(ii)] \emph{uniformly $\hK$-stable} if there exists $\sigma>0$ such that 
        $\mm_\na(\f)\ge\sigma\inf_{c\in\R}\dd_{1,\na}(\f,c)$ 
        for $\f\in\cE^1_\na$.
    \end{itemize}
\end{defi}
Here we view $\R$ as a closed subset of $\cE^1_\na$, consisting of constant potentials. As we discuss in~\S\ref{sec:KandKhat} below, (uniform) $\hK$-stability implies (uniform) K-stability; we do not know if the two are equivalent in general.

%
%
\subsection{A version of the YTD conjecture}
The result below is a slightly more precise version of Theorem~A in the introduction.
\begin{thm}\label{thm:ytd1}
    Assume $\Aut^0(X,L)$ is trivial. Then the following conditions are equivalent:
\begin{itemize}
    \item[(i)] there exists a cscK metric in $c_1(L)$;
    \item[(ii)] $(X,L)$ is $\hK$-stable;
    \item[(iii)] $(X,L)$ is uniformly $\hK$-stable.
\end{itemize}
\end{thm}
Roughly speaking, this says that the existence of a cscK metric is equivalent to the trivial metric on $L$ being a `non-Archimedean cscK metric' (in the sense of minimizing the non-Archimedean Mabuchi functional). 
%
%
\subsection{CscK metrics and the Mabuchi functional}
The proof of the YTD conjecture is variational, and a nontrivial part of it involves understanding the relationship between cscK metrics and the Mabuchi functional on $\cE^1$. Indeed we have the following result, \cf~\cite{DL20,CC2}:
\begin{thm}\label{thm:ytdgeod}
    If $\Aut^0(X,L)$ is trivial, then the following conditions are equivalent:
    \begin{itemize}
        \item[(i)] there exists a (unique) cscK metric in $c_1(L)$;
        \item[(ii)] $(X,L)$ is \emph{geodesically stable}, \ie $\lim_{t\to\infty} t^{-1}\mm(\phi_t)\ge0$ for any psh geodesic ray $\{\phi_t\}$ in $\cE^1$, with equality iff $\phi_t=\phi_0+c t$ for $c\in\R$;
        \item[(iii)] the Mabuchi functional is \emph{coercive} on $\cE^1/\R$, \ie there exist $\sigma,C>0$ such that 
        $\mm(\phi)\ge\sigma\inf_{c\in\R}\dd_1(\phi,\phi_\refe+c)-C$ for $\phi\in\cE^1$.
    \end{itemize}
\end{thm}
We will not attempt to recap the proof of Theorem~\ref{thm:ytdgeod}, but it relies on the completeness and Busemann convexity of $\cE^1$, together with the following key properties of the Mabuchi functional:
\begin{itemize}
    \item $\mm$ is convex along psh geodesics in $\cE^1$~\cite{BerBer,BDL17};
    \item $\mm$ is strongly lsc in the sense that the intersection of any sublevel set of $\mm$ with a closed ball is compact~\cite{BBEGZ};
    \item any minimizer of $\mm$ on $\cE^1$ lies in $\cH$~\cite{BDL20,CC1,CC2};
    \item minimizers of $\mm$ are unique up to an additive constant~\cite{BerBer} (recall that $\Aut^0(X,L)$ is assumed to be trivial).
\end{itemize}
%
%
\subsection{Proof of Theorem~\ref{thm:ytd1}}
That~(iii) implies~(ii) is a consequence of $\R\subset\cE^1_\na$ being a closed subset. To prove the remaining implications, we use Theorem~\ref{thm:ytdgeod} together with the asymptotics of the Mabuchi functional in~\S\ref{sec:entasym}.

Indeed, \eqref{equ:mabasym} immediately implies that $(X,L)$ is $\hK$-stable iff it is geodesically stable. In view of Theorem~\ref{thm:ytdgeod}, we therefore have that $(X,L)$ being $\hK$-stable implies the existence of a cscK metric in $c_1(L)$. 

It remains to prove the implication~(i)$\Rightarrow$(iii), so suppose $c_1(L)$ contains a cscK metric. By Theorem~\ref{thm:ytdgeod} there exist $\sigma,C>0$ such that 
$$
\mm(\phi)\ge\sigma\inf_{c\in\R}\dd_1(\phi,\phi_\refe+c)-C
$$
for $\phi\in\cE^1$. Applying this along a maximal geodesic ray $\{\phi_t\}$ directed by $\f\in\cE^1_\na$ and using~\eqref{equ:mabasym}, one can show
$$
\mm_\na(\f)\ge\sigma\inf_{c\in\R}\dd_{1,\na}(\f,c)
$$
on $\cE^1_\na$, which proves that $(X,L)$ is uniformly $\hK$-stable.
%
%
%
%
\section{Extensions and further remarks}\label{sec:extensions}
Theorem A in the Introduction is a special case of the results in~\cite{wytd}. Let us briefly comment on more general situations.
%
%
\subsection{Generalizations}
Above we assumed that $\Aut^0(X,L)$ is trivial. In general, a cscK metric in $c_1(L)$ is unique up to the action of $\Aut^0(X,L)$~\cite{BerBer}, and its existence is equivalent to \emph{$\hK$-polystability}, and also to \emph{uniform $\hK$-polystability}, in which the set $\R\subset\cE^1_\na$ of constant potentials is replaced by the set $\cP_\R\subset\cE^1_\na$ of \emph{real product test configurations}. We can also look for cscK metrics invariant under a given compact subgroup $S\subset\Aut^0(X,L)$. Finally, we can more generally treat \emph{extremal metrics} in the sense of Calabi, and even \emph{weighted extremal metrics} in the sense of Lahdili~\cite{Lah19}. We refer to~\cite{wytd} for the precise statements.

In another direction, one can consider the `transcendental' case of an arbitrary Kähler class on a compact Kähler manifold. The relevant notion of uniform K-stability in that context was introduced in independent works of Dervan--Ross and Sj\"ostr\"om-Dyrefelt \cite{DR,SD}, while those of Berkovich space and $\hK$-stability where more recently developed in Mesquita-Piccione's PhD work~\cite{MP}, which further adapted C.~Li's approach to show that uniform $\hK$-stability implies the existence of a cscK metric. Here the converse direction remains to be proved, but a valuative/divisorial criterion for $\hK$-stability was obtained in~\cite{MW} (see below). 
%
%
\subsection{K-stability}\label{sec:KandKhat}
The `classical' notions of \emph{K-stability}~\cite{Tia97,Don02} and \emph{uniform K-stability}~\cite{Der16,BHJ1} are formulated in terms of test configurations: they are equivalent to conditions~(i) and~(ii) in Definition~\ref{defi:hKstab}, but with $\cE^1_\na$ and $\R$ replaced by $\cH_\na$ and $\Q$, respectively. Thus they are weaker than $\hK$-stability and uniform $\hK$-polystability, respectively, and one may ask if they are equivalent. A positive answer for the uniform notion would follow if the entropy regularization conjecture were to hold; see~\S\ref{sec:cscK}. When $\Aut^0(X,L)$ is nontrivial, such a regularization result would also prove that uniform $\hK$-polystability is equivalent to \emph{uniform K-polystability}, a suitable quantitative lower bound of $\mm_\na$ on $\cH_\na$.
%
%
\subsection{The approach of Darvas and Zhang}
T.~Darvas and K.~Zhang~\cite{DZ25} independently found a different YTD-type correspondence: in the case when $\Aut^0(X,L)$ is trivial, the existence of a cscK metric in $c_1(L)$ is equivalent to $(X,L)$ being \emph{$\mathrm{K}^\beta$-stable}, a condition on test configurations defined in terms of a regularized Mabuchi functional $\mm^\beta_\na$, depending on a parameter $\beta>0$, with better continuity properties.
%
%
\subsection{Valuative criterion}
Even when $\Aut^0(X,L)$ is trivial, verifying whether a polarized smooth variety $(X,L)$ is $\hK$-stable or uniformly $\hK$-stable appears to be quite intractable, as the set $\cE^1_\na$ of test objects is enormous. 

In the case of Fano manifolds, thanks to outstanding technological progress based in part on the Minimal Model Program (see~\cite{XuBook}), the task of verifying K-(poly)stability is now more feasible, as seen in the classification~\cite{ACC+,Fuj23}. A key ingredient here is the Fujita--Li valuative criterion for K-stability~\cite{FujVal,LiVal}, which was later extended to the general polarized case by Dervan--Legendre~\cite{DL23} and then the authors~\cite{nakstab2}. In this generality, the valuative criterion consists in checking the positivity of the \emph{$\beta$-invariant} of all divisorial measures. However, while computing this invariant for a given measure is made reasonably explicit by a formula of C.~Li and Mesquita-Piccione--Witt Nystr\"om~\cite{LiLect,MP}, the valuative criterion remains impractical at this point to check $\hK$-stability, as the set of divisorial measures is still an infinite one.
%
%
\subsection{Moduli spaces}
One major motivation for the YTD conjecture is the construction of moduli spaces of (suitably stable) polarized varieties, something that has indeed been successfully carried out for K-stable Fano manifolds/varieties in a very impressive series of works (see~\cite{XuBook}). Analytic moduli spaces of compact Kähler manifolds equipped with a cscK metric were constructed by Fujiki--Schumacher~\cite{FS} in the case of trivial automorphism group, and by Dervan--Naumann in general~\cite{DN}. However, fundamental questions towards a more refined understanding of the problem, such as the Zariski openness of (uniform) $\hK$-stability in families, remain widely open at this point.
%
%
%
%

\end{document}